\documentclass[12pt,a4paper]{article}
\usepackage[T2A]{fontenc}
\usepackage{amssymb}
\usepackage{latexsym}
\usepackage{amsmath}
\usepackage[active]{srcltx} 
\usepackage{amsmath,amssymb,amsthm} 

\usepackage[normalem]{ulem}
\usepackage{color}
\usepackage{colordvi,multicol}

\newtheorem{theorem}{Theorem}
\newtheorem{pro}{Property}
\newtheorem{corollary}{Corollary}
\newtheorem{lemma}{Lemma}

\newtheorem{defn}{Definition}

\def\Z{\mathbb{Z}}    
\def\X{\mathbb{X}}

\renewcommand{\le}{\leqslant}

\renewcommand{\ge}{\geqslant}

\newcommand{\oo}[1]{\overline{#1}}
\newcommand{\ww}{\widetilde}

\title{Structural Properties of \\
 Conditioned Random Walks  on  
  Integer Lattice{s} 
\\ with Random Local Constraints\footnote{Research is supported
by RSF research grant No. 17-11-01173} }
\author{Sergey Foss\footnote{Email address: sergueiorfoss25@gmail.com}\\Heriot-Watt University and \\ Novosibirsk State University and\\ Sobolev Institute of Mathematics 
\and Alexander Sakhanenko\footnote{Email address: aisakh@mail.ru} \\ Novosibirsk State University and \\ Sobolev Institute of Mathematics}
\date{\today}

\begin{document}

\maketitle

\begin{center}
{\it \bf To the memory of Vladas}
\end{center}

\begin{abstract}
 We consider a random walk on a multidimensional integer lattice with random bounds on local times, conditioned on the event that it hits a high level before its death.
We introduce an auxiliary ``core'' process that has a regenerative structure and plays a key role in our analysis.
We obtain a number of representations for the distribution of the random walk in terms of the
similar distribution of the ``core'' process.
 Based on that, we prove a number of limiting results by letting the high level to tend to infinity. In particular, we generalise results for a simple symmetric one-dimensional random walk obtained earlier in the paper by Benjamini and Berestycki (2010). 
 We are thankful to Vladas Sidoravicus who has discussed the paper with Sergey Foss.
 \end{abstract}

{\bf Keywords:} Conditioned Random Walk, Bounded Local Times, Regenerative Sequence, Potential Regeneration, Separating Levels, Skip-Free Distributions

\section{Introduction}
Consider a $d$-dimensional random walk
\begin{gather}                                                                               \label{i1}
S_t=(S_t[1],\dots,S_t[d]) = S_0+ \sum_{j=1}^t\xi_j, \quad t=0,1,2,\dots,
\end{gather}
on the integer lattice $\Z^d$, 
where $\xi_j=(\xi_j[1],\ldots,\xi_j[d]),$ $j=1,2,\ldots$ are i.i.d. 
random vectors  
 that do not depend on the initial value $S_0$.
The random variable
\begin{equation}                                                                                                             \label{i3}
L_t(x) = \sum_{j=0}^t {\bf 1}{\{S_j=x\}}, \quad x\in \Z^d, 
\end{equation}
counts the number of visits to (or local time at) state $x$ by time $t=0,1,2,\ldots$.
 We assume that, for each $x\in {\Z}^d$, the number of possible/allowed visits to
 state~$x$ is limited above by a counting number $H(x)\ge 0$.       
   Let
\begin{equation}                                                                                                               \label{i4}
T_*=\inf \{ t\ge0  : L_{t}(S_{t}) >H(S_{t})\}\le\infty
\end{equation}
be the first time when the number of visits to any state exceeds its upper limit.
 If $ T_*$ is finite, we assume that the random walk is ``killed''
 at the time instant~$ T_*$ (or it ``dies'', or ``freezes'' 
at time $ T_*$).

Thus, we consider  a multidimensional integer-valued random walk in a 
changing random environment, where initially 
each point $x$   is characterized by a random number $H(x)$ of allowed  visits to it. 
At any time $t$,  the random walk jumps from $x=S_t$ to $S_{t+1}$ and changes the environment at point $x$ by decreasing the number of remaining allowed visits by 1.
As a natural example, consider a model of a random walk on atoms of a ``harmonic crystal'' (see, e.g., \cite{BDG} and \cite{BB}).
An electron jumps from one atom to another, 
taking from a visited atom for the next jump a fixed unit of energy, which cannot be recovered. 
Thus, if $S_{t}$ is a position of the electron at time $t$, then a unit of energy is sufficient 
for it to have a next jump to position $S_{t+1}=S_{t}+\xi_{t+1}$, which may be in any direction from~$S_{t}$.
We interpret the first coordinate $S_t[1]$ of $S_t$ as its {\it height} and assume further that the height cannot increase by more than one unit:
\begin{gather}                                                                             \label{i5}
\xi_t[1]\le1\quad a.s.
,\qquad t=1,2,\dots.
\end{gather}

When the electron arrives at an atom with insufficient energy level, it ``freezes'' there. 
We may formulate two natural tasks. Firstly, to find the asymptotics, as $n\to\infty$, of the probability of the event $B_n$
that  the electron reaches the level $n$ before it ``freezes'', i.e.
\begin{equation}                                                                     \label{i7}
B_n :=\{\alpha(n)< T_*\} 
\end{equation}
with $\alpha(n)$ being the hitting time of the level $n$:
\begin{align} \label{i8}
\alpha(n) := \inf \{ t\ge0  :  S_t[1]\ge n\} =\inf \{ t\ge0  :  S_t[1]= n\}.
\end{align}  
where the latter equality follows from the skip-free property \eqref{i5}.

Secondly, given that the electron is still active by the time of hitting  level $n$,
 a  question of interest is the asymptotic, as $n$ increases, of the conditional distribution of the electron's sample path.

To clarify the presentation, we will use the low-case ``star'' in the probability ${\mathbf P}_*(\cdot)$ in order to
underline the influence of the random environment.
We omit the ``star'' in ${\mathbf P}(\cdot)$
if the environment is not involved.  

 In \cite{BB}, a simple symmetric one-dimensional random walk on the integers  
 (``one-dimensional atoms'') has been considered  
under the assumptions  that 
\begin{align*}       
{\mathbf P}(\xi_1=1)=1/2={\mathbf P}(\xi_1=-1),
\quad S_0=0 \quad \mbox{and}\quad 
H(x)=L_0 =const \ge 2 
\end{align*}
for all $x\in\Z$.
The latter means that  initially each atom has a fixed (the same for all)
amount of energy~$L_0$. 
The authors showed that 
\begin{align*}    
{\mathbf P}_*(B_n)\sim const\cdot q^n
\quad \mbox{as}\quad n\to \infty,
\quad \mbox{where}\quad 0<q<1.
\end{align*}
Based on that, they proved (see Theorem 5 in \cite{BB}) 
convergence of the conditional distributions:
\begin{equation}                                                                               \label{i10}
{\mathbf P}_* ((S_0,\ldots,S_k)\in A \ | \ B_n )  \to
{\mathbf P} ((\overline{S}_0,\ldots,\overline{S}_k)\in A),
\end{equation}
for any $k=0,1,2,\ldots$ and all $A \subset \Z^{(k+1)\times d}$,
where $\Z^{K\times d}$ denotes the space of vectors $\vec x=(x_1,\dots,\vec x_K)$
having $d$-dimensional vectors as their components.
Further, it was shown in \cite{BB} that the limiting sequence $\{\oo{S}_k\}$ in (\ref{i10}) has a  regenerative structure
(see Definition~\ref{def3} 
below 
 for details)
and increases to infinity with a linear speed, i.e.
\begin{equation}                                                                               \label{i11}
\oo S_n/n\to1/\mu
\quad \mbox{a.s. as}\quad n\to \infty,
\quad \mbox{where}\quad 1\le\mu<\infty.
\end{equation}

In our paper, we consider a multivariate random walk on the integer lattice with random local constraints. We generalise the model of \cite{BB} in three directions: we consider 
more general distributions of jumps, many dimensions,
and  random local constraints.
We develop the approach introduced in \cite{BB}, with a number of essential differences. 
The main difference is that we first focus on the analysis of the structure 
of the initial random walk $\{S_t\}$. In particular, we introduce a notion  of  $n$-separating levels which often exist in our model.
The analysis of properties of such random levels allows us to introduce a sequence 
of random vectors $\{\oo S_t\}$ with specially chosen joint distribution.
We call  $\{\oo S_t\}$ the {\it core} random sequence, or the {\it core} random process.

There are several advantages of studying the core process. We show that its 
structure 
 (a) does not involve any counting constraints, (b) does not involve an environment, (c) operates with
proper distributions only, and 
(d) the core process has a (strongly) regenerative structure
with an infinite sequence of  
random regenerative levels $\{\oo\nu_i\}$
(see Definition~\ref{def3} for details).

We obtain a number of interesting representations for the conditional distribution of the random walk $\{S_t\}$ in random environment $\{H(x)\}$, linked to the distribution of the core sequence $\{\oo{S}_t\}$. 
These representations allow us to obtain a number of novel results. 
For example, we show that
\begin{align}                                                                                                       \label{i12}
{\mathbf P}_* (B_n) =\psi_0q^n {\mathbf P}(\overline{B}_n),
\quad \mbox{where} \quad
\overline{B}_n :=\cup_{m=0}^n\{ \overline\nu_m=n \},
\end{align}
for well-defined positive constants $\psi_0$ and $0<q\le1$, and that 
\begin{equation}                                                                                                   \label{i13}
{\mathbf P}_* ((S_0,\ldots,S_k)\in A \ | \ B_n ) =
{\mathbf P} ((\overline{S}_0, \ldots,\overline{S}_k)\in A \ | \ \oo{B}_n),
\end{equation}
for any $n\ge k=0,1,2,\ldots$   and all $A \in \Z^{(k+1)\times d}$.
Here event $\oo{B}_n$ occurs  iff $n$ coincides with one of the 
regenerative levels of the core random walk.

Finally, we obtain the desired limiting result (\ref{i10}) as a simple corollary 
of (\ref{i13}), which is 
a generalisation of Theorem 5 in \cite{BB}.

We have to mention that a number of known results for conditioned random walks that do not have local-time constraints (see, e.g. \cite{BD} and \cite{A}) may be represented, in some particular cases, as corollaries of our results, see Section 7.2 for detailes.

There is a number of publications on random walks with constraints on local times. 
We have already mentioned papers \cite{BB} and \cite {BDG}. The paper \cite{BB} was, in fact, the initial point of our studies,
and we have made a number of preliminary observations in \cite{SF} where we considered a reasonable one-dimensional 
generalisation of the discrete-time model in \cite{BB} with non-random boundary constraints.  Papers \cite{BB2} and \cite{KS} deal with a different problem: they consider 
a random walk on the line (see also \cite{B} for a generalisation onto a class of Markov processes), assuming that the initial energy level $H(x)$ of a point $x>0$ is a deterministic function of $x$ that increases to infinity with $x$. These papers analyse 
 recurrence/transience properties of the random walk that depend on the shape of the function $H(x)$.
A generalisation of the model onto random trees may be found in \cite{BGMS}.    
Papers \cite{FZ}, \cite{MS} and \cite{Ku} are more distant, they discuss
unconditioned regenerative phenomena that depend on an infinite future, in a number of situations.  

To conclude, in the present paper we  provide a unified treatment of the
conditional regenerative phenomenon in a class of multivariate random walks on the integers with changing random constraints on the numbers of visits. 

The paper is organized as follows. In Section 2, we introduce 
the main assumptions on the model and the notions of separating and regenerative 
levels.
In Section 3, we first introduce and discuss the 
structure of
the model connected with the existence of random $n$-separating levels.
After that 
we  describe the core random sequence and its structure, and, finally,  formulate the Representation Theorem and  
limiting results as its Corollaries. Then Sections 4 -- 6 are devoted to the proofs. 
We have to note that, in the proof of the main auxiliary result, the Key Theorem, we follow
the approach developed in \cite{BB}. We conclude with Section 7 containing a few remarks.


\section{Main Assumptions and Definitions}

\subsection{Basic Assumptions}

For $n\in\Z$, introduce a half-space of $\Z^d$ 
\begin{align}                                                                              \label{a0}
\Z^d_{n+}:=\{x=(x[1],x[2],\dots,x[d])\in \Z^d:x[1]\ge n\}.
\end{align}
The following assumptions $(A1)-(A3)$,  are supposed to hold throughout the
paper.

${\bf (A1)}$.   
The increments $\{\xi_t:  t\ge 1\}$ of the random walk $\{S_t \}$ from \eqref{i1} are i.i.d. random vectors taking values in $\Z^d$, and their first components  
have a skip-free distribution: 
\begin{equation*} 
\sum_{k=-\infty}^1 {\mathbf P} (\xi_1[1]=k) = 1 
\quad \mbox{and}\quad
{\mathbf P} (\xi_1[1]=1)>0.
\end{equation*}
 
 ${\bf (A2)}$. 
The random constraints $\{H(x),\, x\in \Z^d\}$ 
are non-negative integer-valued random variables which
may take the infinite value: for any $x\in\Z^d$, 
\begin{equation}                                                                                                 \label{a2}
\sum_{l=0}^\infty{\mathbf P}(H(x)=l)+{\mathbf P}(H(x)=\infty)=1.
\end{equation}
Moreover, the next three families of random variables
\begin{equation*}    
\{S_0;\ H(x), \,  x\notin\Z^d_{0+}\}, \quad \{\xi_i,\ i\ge 1\} 
\quad \mbox{and}\quad
\{H(x), \  x\in\Z^d_{0+}\}
\end{equation*}
are mutually independent, \ $S_0[1]\le0$ \ a.s. \ and \ ${\mathbf P}_* (B_0)>0$.

${\bf (A3)}$. 
The family  $\{H(x):x\in\Z^d_{0+}\}$ consists of i.i.d. random variables with 
\begin{equation*}   
{\mathbf P}(1\le H(0)\le\infty)=1.
\end{equation*}

We may interpret Assumption (A3) as follows: at time $t=0$ the environment in~$\Z^d_{0+}$ is stochastically homogeneous, so is ``virgin'' (see, also, Remark 7.1). Then condition ${\mathbf P}_*(B_0)>0$ in Assumption (A2) may be read as ``the random walk $S_t$ arrives at the 
virgin domain of the random environment with a positive probability.''

Assumptions $(A1)$ -- $(A3)$ yield that, for any $n\ge 0$, 
\begin{gather}                                                                                                                \nonumber
 {\mathbf P}_* (B_{n})\ge{\mathbf P}_* ( \alpha(0)<T_*,\,\xi_{\alpha(n)+j}[1]=1,\,H(S_{\alpha(n)+j})>0,\,j=1,\ldots,n)
 \\                                                                                                               \label{a3}
 = {\mathbf P}_* (B_{0}) {\mathbf P}^{n} (\xi_1[1]=1)>0,
\end{gather}
where the events $B_n$ were introduced in (\ref{i7}). Thus, for all $n\ge0$ the event $B_n$
occurs with positive probability and hence, as we can see later, all conditional probabilities in all our main assertions are well defined.

\subsection{Technical Assumption and Comments}
We have certain flexibility in the initial value $S_0$ and in the random environment $\{H(x)\}$ outside the set $\Z^d_{0+}$. 
Recall that we use notation ${\mathbf P}_*(\cdot)$ for probabilities of events where the environment is involved.
We will also use special notation, ${\mathbf P}_0$ and ${\mathbf P}_+$, for two particular environments
{when $S_0=0$}. 
  For any event $B$, let
\begin{gather}                                                                                      \label{a4}
{\mathbf P}_0(B):={\mathbf P}_*(B\,|\, S_0=0,\ H(y)=0\ \ \forall  y\notin\Z^d_{0+}\},
\\                                                                                      \label{a4+}
{\mathbf P}_+(B):={\mathbf P}_*(B\,|\, S_0=0,\ H(y)=\infty\ \ \forall  y\notin\Z^d_{0+}\}.
\end{gather}
In (\ref{a4}), {\it it is prohibited for the random walk to visit any states} $ y\notin\Z^d_{0+}$, and 
(\ref{a4+}) corresponds to the case where there is 
 {\it no restrictions on the number of visits to any of the states} $ y\notin\Z^d_{0+}$.
Clearly, 
\begin{gather}                                                                                                   \label{a5}
{\mathbf P}_+(B_0)={\mathbf P}_0(B_0)=1\quad \mbox{and}\quad
{\mathbf P}_+ (B_{n})\ge{\mathbf P}_0 (B_{n})\ge {\mathbf P}^{n} (\xi_1[1]=1)>0\ \ \forall n\ge0.
\end{gather}

For the classical random walk (no environment), introduce two stopping times:
\begin{equation*}                
\beta_0:=\inf \{ t>0  : \xi_1[1]+\dots+\xi_t[1]=0\}\le\beta_{0,0}:=\inf \{ t>0  : \xi_1+\dots+\xi_t=0\}\le\infty.
\end{equation*}
We will need the following assumption:

 ${\bf (A4)}$ If ${\mathbf P}(\beta_0<\infty)=1$ then ${\mathbf P}(\beta_{0,0}<\infty)>0$.
And
\begin{equation*}                                                                                                   \label{a7}
\mbox{if} \   {\mathbf P}(\beta_0=\beta_{0,0}<\infty)=1, \ \mbox{then} \ 
{\mathbf E}H(0)<\infty.
\end{equation*}

It is clear that assumption (A4) is fulfilled in the following cases:
\\
(a) \ ${\mathbf E} \xi_1[1]\ne 0$ (including the cases ${\mathbf E} \xi_1[1]> 0$  and \ $-\infty \le {\mathbf E} \xi_1[1] < 0$);
\\
(b) \ ${\mathbf E} \xi_1[1]= 0$ and $0<{\mathbf P}(\beta_{0,0}<\infty)<1$;
\\
(c) \ ${\mathbf E} \xi_1[1]= 0$, \ ${\mathbf P}(\beta_{0,0}<\infty)=1$
and ${\mathbf E}H(0)<\infty$.

Thus, our results do not work only in the next two cases:
\\
(d) \ ${\mathbf E} \xi_1[1]= 0$ and ${\mathbf P}(\beta_{0,0}<\infty)=0$;
\\
(e) \ ${\mathbf E} \xi_1[1]= 0$, \ ${\mathbf P}(\beta_{0,0}<\infty)=1$, \ and  ${\mathbf E}H(0)=\infty$.

Note the case (d) is degenerate in the spirit of our paper,
since it corresponds to the situation where the random walk visits each state at most once.

Note also  that the cases (c) and (e) 
relate to essentially one- or two-dimensional random walks only.  

\subsection{Separating and Regenerative Levels}
For a finite or  infinite sequence  
$\vec y = (y_0, y_1, y_2, \ldots)$  of $\Z^d$-valued vectors and for  any $n\ge0$,  we let
\begin{gather}                                                                                                                     \label{a8}
\alpha(n|\vec y):=\inf\{t\ge0: y_t[1]\ge n\}\le\infty, 
\end{gather}  
where
$y_t[1]$ is the first coordinate of $y_t$, for  $t=0,1,\ldots.$
Here and throughout the paper, we follow the standard conventions that 
\begin{equation}                                                                                                          \label{a9}
 \inf\emptyset=\infty,  \quad  \sup\emptyset=-\infty
  \quad \mbox{ and}\quad\sum_{k\in\emptyset}a_k=0. 
\end{equation}

\begin{defn}                                                              \label{def1}
A number $k\ge0$  is a ``separating level'' of the  sequence $\vec y$ if 
\begin{gather*}                                                                                                                    
\alpha (k|\vec y)<\infty \quad \mbox{and} \quad 
\max_{0\le t<\alpha(k|\vec y)}y_t[1]<k=y_{\alpha(k|\vec y)}[1]\le\inf_{t>\alpha(k|\vec y)}y_t[1].
\end{gather*}
\end{defn}

\begin{defn}                                                              \label{def2}
A number $k\in\{0,1,\ldots,n\}$  is an ``$n$-separating level''  of the  sequence $\vec y$ if
\begin{gather*}                                                                                                                    
\sup_{0\le t<\alpha(k|\vec y)}y_t[1]<k=y_{\alpha(k|\vec y)}[1]\le\inf_{\alpha(k|\vec y)<t<\alpha(n|\vec y)}y_t[1]
\quad \mbox{and} \quad \alpha(n|\vec y)<\infty.
\end{gather*}
\end{defn}

For $n\ge 0$, let $\eta (n|\vec y)+1$ counts the number of $n$-separating levels;
and let $\varkappa(n|\vec y)$ be the supremum of all $k<n$ such that
$k$ is an $n$-separating level. 

These levels play an important role in our analysis. One can see that 
if $k$ is an $n$-separating level, then it may not be an $N$-separating level
for $N>n$ and, hence, it may be not a separating level. For example, $k=n$ is always the last $n$-separating level if
$\alpha(n|\vec y)$ is finite, 
but it is not an $(n+1)$-separating level if $y_{\alpha(n)+1}[1]< 0$.

In what follows, a ``block'' is  any collection of random variables
that may contain a random number of these variables.  
\begin{defn}    \label{def3}
A random sequence $\oo{S}=(\oo{S}_0,\oo{S}_1,\dots)$ is {\it 
 strongly regenerative with regenerative levels}
 $\overline{\nu}_0<\overline{\nu}_1< \ldots < \overline{\nu}_n < \ldots$,
 if $\{\oo{\nu}_i\}$ is an infinite sequence of proper integer-valued random variables such that, firstly,
 the following ``blocks'' of random variables 
  \begin{gather*}   
\{\oo{\nu}_i-\oo{\nu}_{i-1},\, \oo{\alpha}(\oo{\nu}_i)-\oo{\alpha}(\oo{\nu}_{i-1}),\,
(\oo{S}_{\oo{\alpha}(\oo{\nu}_{i-1})+t}-\oo{S}_{\oo{\alpha}(\oo{\nu}_{i-1})},\, 
 t=1,2,\ldots, \oo{\alpha}(\oo{\nu}_i)-\oo{\alpha}(\oo{\nu}_{i-1}))\},\ \  i\ge1,
\end{gather*}
 are i.i.d. and do not depend on the initial ``block'' \ \
$
\{\oo{\nu}_0,\oo{\alpha}(\oo{\nu}_0),(\oo{S}_t; t\le \oo{\alpha}(\oo{\nu}_0))\},
$
and, secondly,
\begin{gather*}       
  \inf_{t\ge \oo{\alpha}(\oo\nu_i)} \overline{S}_t[1]=\overline{S}_{\oo{\alpha}(\oo\nu_i})[1]=\oo\nu_i
  >\sup_{ 0\le t< \oo{\alpha}(\oo\nu_i)} \overline{S}_t[1],\quad i=0,1,2,\dots.
   \end{gather*}
\end{defn}

We then say that $\oo{\alpha}(\oo{\nu}_i)$ is the {\it regenerative time} that corresponds to regenerative level~$\nu_i$. One can view 
$n$-separating levels as ``potential candidates'' for regenerative levels
and talk about ``potential regeneration''.

\section{Main Results}
In Subsection 3.1 we introduce a renewal equation for the random walk with local constraints and introduce its splitting into random blocks. In Subsection  3.2 we present the Key Theorem  and introduce a sequence of independent blocks related to the core sequence. Based on that, we provide a formal definition of the core process in Subsection~3.3 .
After that we present our main results in Subsections 3.4 and 3.5.

\subsection{
On the structure of the random walk}
Note that earlier notation \eqref{i8}  matches \eqref{a8}  as follows: $\alpha(n)=\alpha(n| S)$,
for ${S}=({S}_0,{S}_1,\dots)$. For each $n\ge0$, we let 
\begin{gather}                                                                                                       \label{a15}
 \eta_*(n):=
\begin{cases}
\eta(n|S)
 \ \mbox{if}\ \ \alpha(n)<T_*(n),
\\
-1, \qquad\quad \mbox{otherwise},
\end{cases}
   \mbox{and}\ \
 \varkappa_*(n):=
\begin{cases}
\varkappa(n|S), \  \mbox{if}\ \ \eta_*(n)\ge1,
\\
-\infty, \ \ \ \ \mbox{if}\quad \eta_*(n)<1.
\end{cases}
\end{gather}

So,  $\eta_* (n)+1$ counts the number of $n$-separating levels in the case where the event $B_n = \{\eta_*(n)\ge0\}$ occurs.  
Note that if the event $B_n$ occurs, then $k=n$ is the largest  $n$-separating level, and $k=\varkappa_*(n)$ is the second largest $n$-separating level, if it exists, i.e. when $\eta_* (n)\ge1$.
Clearly, 
\begin{gather}                                                                       \label{a17}
\{\varkappa_*(n)>-\infty\}=\{0\le\varkappa_*(n)\le n-1\}=\{\eta_*(n)\ge1\}\subset \{\eta_*(n)\ge0\}=B_n.
\end{gather} 
Further, ${\mathbf P}_0(\eta_*(n)=0)=1$ 
because, under the ``$0$-environment'', level $0$ is $n$-separating for any $n$ such that $\alpha(n)<T_*$. 

 The random walk under consideration has the following renewal-type  Property. 
\begin{pro}                                                                                                          \label{pro1}
Under the assumptions $(A1)$ -- $(A3)$, for any $n>k\ge0$,
\begin{gather*}                                                                                                          \label{a19}
{\mathbf P}_* (B_n, \varkappa_*(n)=k)={\mathbf P}_*(B_k) \cdot
{\mathbf P}_{0} (\varkappa_*(n-k)=0),
\end{gather*}
and then  the following renewal equation holds:
\begin{gather}                                                                                                          \label{a20}
{\mathbf P}_* (B_n)={\mathbf P}_* (\eta_*(n)=0)+\sum_{k=0}^{n-1} {\mathbf P}_*(B_k) \cdot
{\mathbf P}_{0} (\varkappa_*(n-k)=0),\quad n=1,2,\dots.
\end{gather}
In particular,
\begin{gather}                                                                                                          \label{a21}
{\mathbf P}_0 (B_n)=\sum_{k=0}^{n-1} {\mathbf P}_0(B_k) \cdot
{\mathbf P}_{0} (\varkappa_*(n-k)=0),\quad n=1,2,\ldots.
\end{gather}
\end{pro}

For $n> 0$ with $\eta_*(n)\ge0$, let
\begin{gather*}                                                                                                    \label{a22}
0\le\nu_0(n)<\ldots<\nu_{\eta_*(n)}(n)=n 
\end{gather*}
be the sequence of all $n$-separating levels (where 
$\nu_0(n)=\nu_{\eta_*(n)}(n)=n$ if $\eta_*(n)=0$). In the case $\eta_*(n)\ge1$, we may find 
 all $n$-separating levels by the backward recursion:
 \begin{gather*}                                                                                                    \label{a23}
\varkappa_*(\nu_i(n))=\nu_{i-1}(n), \quad i=\eta_*(n),\eta_*(n)-1,\dots,1.
\end{gather*}
For $n>0$ with $\eta_*(n)\ge i\ge1$, we let
\begin{gather*}                                                                                                    \label{a24}
\lambda_i(n):=\nu_i(n)-\nu_{i-1}(n),\quad
T_i(n):=\alpha(\nu_i(n)),\quad
\tau_i(n):=T_i(n)-T_{i-1}(n).
\end{gather*}

We need  more notation.
Introduce the random vectors 
\begin{gather}                                                                                                          \label{a25}
\vec S_K=(S_0,\dots,S_K),\quad \vec S_{K,N}=(S_{K,K+1},\dots,S_{K,N}),\quad N>K\ge0,
\end{gather}
where 
\begin{gather}                                                                          
S_{K,K+t}:=S_{K+t}-S_K=\sum_{j=1}^t\xi_{K+j},\quad t=0,1,\dots.
\end{gather}

On the event $B_n = \{\eta_*(n)\ge 0\}$, introduce 
a random block 
\begin{gather}                                                                                                    \label{a27}
(\nu_0(n),\,T_0(n),\, \vec S_{T_0(n)}).
\end{gather}
This is the {\it initial block} of our random walk. 
Further, if $\eta_*(n)\ge1$, then we may introduce consecutive blocks of  random variables:
\begin{gather}                                                                                                    \label{a28}
(\lambda_i(n), \tau_i(n), 
\, \vec S_{T_{i-1}(n),T_{i}(n)}),\quad i=1,2,\dots, \eta_*(n),
\end{gather}
where $\lambda_i(n)$ is the {\it height} of the $i$-th block and $\tau_i(n)$ its {\it duration}.
Property 1 shows that 
there is a certain 
conditional independence 
of each block in (\ref{a28}) from the previous blocks. 
We present these properties in full in Theorem \ref {T-repr} below.
After that a representation for  the joint distributions of random blocks from (\ref{a27}) 
and (\ref{a28}) will be given in Corollary~\ref{C4}.

\subsection{Key Theorem (the main auxiliary result)}
 The following technical result plays a central role in our studies. 
 It will be proved in Section 5. 

\begin{theorem}                                                                                   \label{T1}
 Under the assumptions $(A1)$ -- $(A4)$,  there exists a 
 number $q\in(0,1]$ such that, 
 \begin{gather}                                                                                \label{q1}
 \sum_{k=1}^{\infty} {\mathbf P}_{0}(\varkappa_*(k)=0)/q^k=1 ,
 \\                                                                                                  \label{q2}
1\le \mu  := \sum_{k=1}^{\infty} k{\mathbf P}_{0}(\varkappa_*(k)=0)/q^k<\infty ,
 \\                                                                                              \label{q3}
0< \psi_0 := \sum_{m=0}^{\infty} {\mathbf P}_* (\eta_*(m)=0)/q^{m} <\infty .
 \end{gather}
\end{theorem}

Properties (\ref{q1}) -- \eqref{q3} allow us to 
introduce 
an infinite sequence  
\begin{gather}                                                                                                  \label{b1}
( \oo{\nu}_0, \oo T_0,  \widetilde S_{\oo T_0})\quad \mbox{and}\quad 
(\oo{\lambda}_i, \oo \tau_i,  \widetilde {Y}_{i,\oo \tau_i}),
\quad i=1,2,\ldots,
\end{gather}
of mutually independent random blocks with special distributions, where 
\begin{gather}                                                                                                  \label{b2}
{\widetilde S}_{\oo T_0}=(\oo S_0,\dots, \oo S_{\oo T_0})
\quad \mbox{and}\quad \widetilde {Y}_{i,\oo \tau_i}=(\oo {Y}_{i,1},\dots, \oo {Y}_{i,\oo \tau_i})
\end{gather}
 are random vectors of random lengths. We determine their distributions step by step.  First, we let
\begin{gather}                                                                                \label{b3}
{\mathbf P} (\oo\nu_0 = k)= {\mathbf P}_* (\eta_*(k)=0)/(\psi_0q^{k}), \quad k=0,1,\dots,
  \\                                                                                              \label{b4}
  {\mathbf P} (\oo{\lambda}_i = l)={\mathbf P}_{0}(\varkappa_*(l)=0)/q^l ,\quad l=1,2,\dots. 
 \end{gather}
 Thus, we have determined the distributions  of random vectors $\oo\nu_0$ and $\oo{\lambda}_i$  as  Cram{\'e}r-type transforms of the characteristics of the initial random walk~$\{ S_t\}$. 
By Theorem~1,  the random vectors $\oo\nu_0$ and $\oo{\lambda}_i$ have proper distributions and
\begin{gather}                                                                          \label{b5}
1\le\mu= {\mathbf E}\oo{\lambda}_1<\infty,\quad 
{\mathbf P} (\oo{\lambda}_1 =1)
={\mathbf P}_{0} (\varkappa_*(1)=0)/q\ge{\mathbf P} (\xi_1[1]=1)/q>0.
\end{gather}
 
 We determine next the distributions of other components of the vectors in  (\ref{b1}). 
We let 
\begin{gather}                                                                                                              \label{b6}
 {\mathbf P} ( \oo T_0=K, {\widetilde S}_{K}=\vec y_K | \oo\nu_0 = k) :=
 {\mathbf P}_* ( \alpha(k)=K<T_*,\vec S_K=\vec y_K|\eta_*(k)=0),
\end{gather}
for any $K\ge k+1\ge1$ and $\vec y_K\in\Z^{(K+1)\times d}$; and then 
\begin{gather}                                                                                                              \label{b7}
 {\mathbf P} ( \oo \tau_i=L, {\widetilde {Y}}_{i,L}=\vec x_L|\oo{\lambda}_i = l) :=
 {\mathbf P}_{0} ( \alpha(l)=L<T_*,\vec S_{0,L}=\vec x_L|\varkappa_*(l)=0),
\end{gather} 
for any $L\ge l\ge1$ and $\vec x_L\in\Z^{L\times d}$.

Thus, we have introduced all joint distributions of random elements 
from (\ref{b1}).
All these distributions are proper, since they are determined by proper 
distributions from (\ref{b3}), (\ref{b4}), (\ref{b6}) and (\ref{b7}). 
By  the construction, with probability~1 
\begin{gather*}                                                                                                          \label{b8}
\oo{\nu}_0\ge 0,\ \ \oo{T}_0\ge 0, \ \ \oo{\lambda}_i \ge 1, \ \ \oo{\tau}_i\ge 1,\ \ 
\oo Z_{i,j}\ge0\ \ \forall\ i,j\ge1.  
\end{gather*} 
Moreover, the  random vectors  
 $\{(\oo{\lambda}_i, \oo \tau_i,  \vec Z_{i,\oo \tau_i}),\ \ i=1,2,\ldots\}$
are i.i.d.

\subsection{Sample-path construction of the core random sequence }
Using mutually independent random blocks introduced in (\ref{b1}), we may define random variables
\begin{gather*}                                                                                                  \label{b11}
\oo{\nu}_m=\oo{\nu}_0+\sum_{i=1}^m\oo{\lambda}_i> \oo\nu_{m-1}
,\quad
\oo{T}_m=\oo{T}_0+\sum_{i=1}^m\oo\tau_i> \oo T_{m-1},
\quad m=1, 2,\ldots.
\end{gather*}
Now we introduce random vectors $\oo S_j$ for all $j\ge0$ using an  induction argument. 
For $j\le\oo{T}_0$ they are given in~(\ref{b2}). 
Suppose we have defined $\oo S_j$ for all $j\le\oo{T}_{i-1}$. 
Then we let
\begin{gather}                                                                                                \label{b12}
\oo S_{\oo{T}_{i-1}+j}:=\oo S_{\oo{T}_{i-1}}+\oo {Y}_{i,j}, 
 \quad j=1,\ldots,\oo\tau_i=\oo{T}_{i}-\oo{T}_{i-1}.
\end{gather}
Thus, we have defined  $\oo S_j$ for all $j\le\oo{T}_{i}$. 
Repeating this procedure for all $i=1,2,\dots$
we define random vectors $\oo S_j$ for all $j\ge0$.

Similar to (\ref{a25}), we introduce vectors with multivariate components:
 \begin{gather}                                                                                                \label{b13}
\widetilde  S_{N}=(\oo S_{0},\dots, \oo S_N),\quad
\widetilde  S_{K,N}=(\oo S_{K+1}-\oo S_K,\dots, \oo S_N-\oo S_K),\quad N>K\ge0.
\end{gather}
Consider now the random blocks
\begin{gather}                                                                                                  \label{b14}
( \oo{\nu}_0, \oo T_0,  \widetilde S_{\oo T_0})\quad \mbox{and}\quad 
(\oo{\lambda}_i, \oo \tau_i,  \widetilde S_{\oo T_{i-1},\oo T_i}),
\quad i=1,2,\ldots,
\end{gather}
and note that, by (\ref{b12}) the $i$-th block in (\ref{b14}) coincides
with the $i$-th block in (\ref{b1}).
Thus, all blocks in (\ref{b14}) are mutually independent 
and all of them, but the initial,  are i.i.d.

\subsection{Representation Theorem}

We are now ready to present our main results.
The following statement summarises the main structural properties of the core random sequence 
and provides an inverse formulae for the distributions of the random walk in terms of the core
 process. 

Let $\Z^{ d}_*:=\cup_{n=1}^\infty\Z^{n\times d}$. We consider $\Z^{ d}_*$
as the state space for random sequences of random lengths.

\begin{theorem}                                                                                                        \label{T-repr}
 Under the assumptions $(A1)-(A4)$,  for any set ${\cal A}\subset \Z^d_*$ and for any $n\ge m \ge 0$, 
\begin{gather}                                                                                                               \label{b21}
{\mathbf P}_* (\alpha(n)<T_*,\,  \eta_*(n)=m,\, (S_0,\ldots,S_{\alpha(n)})\in {\cal A})
\\                                                                           \nonumber
=\psi_0q^n{\mathbf P} (\oo\nu(m)=n,\,(\overline{S}_0,\ldots,\overline{S}_{\oo{\alpha}(n)})\in {\cal A}).
\end{gather}
\end{theorem}

Thus, the distribution of the trajectory of the core random sequence has the same support with the distribution of 
the trajectory of the initial  random walk (any finite sample path  has positive probabilities to occur simultaneously  for the core sequence and for the random walk, however these probabilities may differ). 
In particular,
for all $j=1, 2,\ldots$ the following inequalities hold with probability 1: 
  \begin{gather}                                                                                                              \label{b18}
 \oo\xi_j[1] \le1\quad \mbox{and} \quad \oo S_j[1]\le j,  \quad \mbox{where} \quad \oo\xi_j=\oo S_j-\oo S_{j-1} .
\end{gather}

Since $B_n=\{\eta_*(n)\ge0\}$, 
we have from (\ref{b21}) that,
 for any set $A\subset\Z^{(k+1)\times d}$ , 
\begin{gather}                                                                                                               \label{b23}
{\mathbf P}_* ((S_0,\ldots,S_k)\in A,\, B_n)
=\sum_{m=0}^n{\mathbf P}_* ((S_0,\ldots,S_k)\in A,\, \eta_*(n)=m )
\\   \nonumber
=\sum_{m=0}^n\psi_0q^n{\mathbf P} ((\overline{S}_0,\ldots,\overline{S}_k)\in A,\, \overline\nu(m)=n )
=\psi_0q^n{\mathbf P} ( (\overline{S}_0,\ldots,\overline{S}_k)\in A,\, \overline{B}_n )
\end{gather}
Now \eqref{i12} follows from \eqref{b23} with $A=\Z^{(k+1)\times d}$.
Equating the ratio of the left-hand sides of \eqref{b23} and \eqref{i12}  to the ratio of the right-hand sides leads to \eqref{i13}.

Thus, we have obtained 
\begin{corollary}                                                                                                            \label{C1}
Under the assumptions $(A1)$ -- $(A4)$,  equalities (\ref{i12}), (\ref{i13}) and (\ref{b23}) hold. 
\end{corollary}

\subsection{Limiting Results }
Representation \eqref{b23}  allows us to obtain a number of limiting theorems using the standard renewal arguments.
First of all, we can see from  \eqref{b23}   
that  
\begin{gather*}                                                                        \nonumber
{\mathbf P} (\oo{B}(n)\, |\, \oo{\nu}_0=0)=
{\bf P}(\overline\nu_m=n \quad \mbox{for\ some} \quad m\ge0\, |\, \oo{\nu}_0=0)
\\                                                                          \label{b31}
=\sum_{m=0}^n{\bf P}(\oo{\nu}_i=n\, |\, \oo{\nu}_0=0)
=V_n:={\bf I}\{n=0\}+\sum_{m=1}^n{\bf P}\left(\sum\nolimits_{i=1}^m\overline{\lambda}_i=n\right)
\end{gather*}
 is the renewal function of the undelayed renewal process with i.i.d. increments
$\{\overline{\lambda}_i,\,i=1,2,\dots\}$ satisfying (\ref{b5}). 

Now consider the probabilities 
\begin{gather*}                                                                                                              
U_{n}:={\mathbf P} (\oo A_k \cap\oo B_n),
\quad \mbox{where} \quad \oo A_k:=\{(\oo S_0,\ldots,\oo S_k)\in A\},
 \ A\in\Z^{(k+1)\times d}.
\end{gather*}
Note that $S_j[1]\le k\le\oo\nu_k<\oo\nu_i$ for all $0\le j \le k < i$
{ by (\ref{b18})} . Hence, the event $\oo A_k$ does not depend on the  
random variables $\{\overline{\lambda}_i=\oo\nu_i-\oo\nu_{i-1}:i>k\}$. 
Then 
\begin{gather*}                                                                                                              \nonumber 
U_{n,l}:={\mathbf P} \big(\oo B_n\,|\,\oo A_k,\, \oo{\nu}_k=l\le n\big)
\\                                                                      \label{b33}
={\bf P}\Big(n=\overline\nu_m=\sum\nolimits_{i=1}^m\overline{\lambda}_i 
\ \ \mbox{for\ some} \ m\ge0 
\ \Big|\,\oo{A}_k,\, l=\oo{\nu}_k=\sum\nolimits_{i=1}^k\overline{\lambda}_i\Big)
\\                                                                      \nonumber
={\bf P}\Big(n-l=\sum\nolimits_{i=k+1}^m\overline{\lambda}_i \quad \mbox{for\ some} \quad m\ge k\Big)=V_{n-l}.
\end{gather*}
Hence, by the total probability formula,
\begin{equation*}                                                                                                               \label{b34}
U_{n}-{\mathbf P} (\oo A_k \cap \oo B_n,\oo{\nu}_k>n)=\sum_{l=k}^n{\mathbf P} (\oo A_k,\, \oo{\nu}_k=l)\cdot U_{n,l}
=\sum_{l=k}^n{\mathbf P} (\oo A_k,\, \oo{\nu}_k=l)V_{n-l}.
\end{equation*}
Thus, the differences $U_{n}-{\mathbf P} (\oo A_k,\, \oo B_n,\oo{\nu}_k>n)$ satisfy the renewal equation, where
${\mathbf P} (\oo{\nu}_k>n)\to0$ as $n\to\infty$.
So,  
by \eqref{b5}
 and the local renewal theorem, as $n\to\infty$,  
\begin{equation}                                                                                                               \label{b35}
V_n\to1/\mu \ \  \mbox{and} \  \  U_n\to\sum_{l=k}^\infty{\mathbf P} (\oo A_k,\, \oo{\nu}_k=l)/\mu
={\mathbf P} (\oo A_k)/\mu.
\end{equation}
Substituting  \eqref{b35} into  \eqref{b23} and  \eqref{i12}
leads to the following statement 
\begin{gather}                                                                                                               \label{b36}
{\mathbf P}_* ((S_0,\ldots,S_k)\in A,\ B_n)/q^n\to
\psi_0{\mathbf P} ( (\overline{S}_0,\ldots,\overline{S}_k)\in A )/\mu,
\\                                                                                                              \label{b37}
{\mathbf P} (\oo B_n )\to1/\mu
 \quad \mbox{and} \quad 
{\mathbf P}_* ( B_n )/q^n\to\psi_0/\mu.
\end{gather}
In particular, \eqref{i10} takes place. 
Thus,  we have proved the following result. 
\begin{theorem}                                                                                                            \label{C3}
Under the assumptions $(A1)$ -- $(A4)$,  for all $k\ge0$ and any $A \subset \Z^{(k+1)\times d}$ 
convergences  (\ref{b36}), (\ref{b37}) and (\ref{i10}) hold. In addition, by the Strong Law of Large Numbers,
convergence (\ref{i11}) takes place where the number $\mu$ is given by formula (\ref{q2}).
\end{theorem}
We would like to say that Theorem \ref{C3} was the initial aim of our studies. A simple proof of Theorem \ref{C3} (given above) shows the power of Theorem~\ref{T-repr}. 
In \cite{BB}, direct 
analytical arguments have been used to establish for a simple symmetric random walk a limiting result similar to Theorem \ref{C3}.


\section{Proofs of Property 1 and Auxiliary Lemmas }
\subsection{Additional Notation}
In the proofs we will frequently use notation
\begin{gather*}                                                                                                         \label{d6}
H_{t}(x):=H(x)-L_t(x)=H_{t-1}(x)-{\bf 1}\{S_t=x\},\quad t=0,1,2,\dots,
\end{gather*} 
with $H_{-1}(x):=H(x)$. Thus, $H_t(x)$ is the number of allowed  visits to
 state~$x$ after time~$t+0$.
 
We need a number of further notation. Let
\begin{gather*}                                                                          \label{d1}
h(n):=\min_{0\le t\le n}H_{t}(S_t)=\min\{h(n-1),\, H_{n}(S_n)-1\},\quad n=0,1,2,\dots.
\end{gather*}
{ It follows from \eqref{i4} } that, for all $n,N\ge 0$,
\begin{gather}                                                                          \label{d2}
\{T_*>N\}=\{ h(N)\ge0\},\quad B_n=\{\alpha(n)<T_*\}
=\{ h(\alpha(n))\ge0\}=\{ h(\alpha(n)-1)\ge0\}.
\end{gather}
The latter equality follows from condition $H(x)\ge 1$ for $x[1]\ge 0$.

In what follows, we  consider a random walk that starts at
 time $t\ge 0$ from a state $x$, rather that at time $t=0$ from the state $S_0$.  The following notation will be helpful:
\begin{gather}                                                                                \label{d3}
\alpha_t(l) = \inf \{ j\ge 0  : S_{t,t+j}[1]= l\},\quad
h_t(l,x):=\inf_{0\le j< \alpha_t(l)}H_{t+j}(x+S_{t,t+j}),
\\                                                                          \label{D3}
s(t,L):=\inf_{0\le j< L}S_{t,t+j}[1],
\quad s_t(l):=s(t,t+\alpha_t(l)),
\end{gather}
for $t,l,L\ge0$, where notation $S_{t,j}:=S_{j}-S_t$ for $t\ge j$ was introduced earlier. Note that $\alpha_0(l)=\alpha(l)$ for all $l\ge0$.

Later on we will use the following properties of notation from (\ref{d3}) and (\ref{D3}):
\begin{gather}                   \nonumber 
\alpha(l+m)=\alpha (l)+\alpha_{\alpha(l)}(m),\quad 
\{s(0, T+\alpha_T(l))\ge0\}=\{s(0,T)\ge0,\,S_T+s_T(l)\ge0\},
\\                                                                                            \label{D5}
\{0<T\le T+\alpha_T(l)<T_*\}=\{T>0,\,h(T-1)\ge0,\,h_T(l,S_T)\ge0\},
\end{gather}
for any random or non-random $T\ge0$ and each $l\ge0$ and $m>0$.

Note that, given $\alpha(n)<\infty$, 
\begin{gather*}                                                                                \label{D7}
S_{\alpha(n)}\in \Z_n^d:=\{x=(x[1],x[2],\dots,x[d])\in \Z^d:x[1]= n\}.
\end{gather*}

\subsection{Shifts of Virgin Environment}

For any $j\ge t\ge0$, introduce  random variable
\begin{equation*}                                                                                                             \label{D11}
L_{t,t+j}(x) = \sum_{k=0}^j {\bf 1}{\{S_{t,t+k}=x\}}, \quad x\in \Z^d, 
\end{equation*}
which, similarly to \eqref{i3}, counts the number of visits to state $x$ within  time interval $(t,t+j]$.
For each $k\ge0$, introduce the following (possibly, improper) random variables:
\begin{gather}                                                                                                       \label{D12}
 H^{(k)}(y)=
\begin{cases}
H(y), \  y\in\Z^d_{k+},
\\
\infty,\, \ \quad y\notin\Z^d_{x[1]+},
\end{cases}
   \mbox{so\ that}\ \
 H^{(x[1])}(x+y)=
\begin{cases}
H(x+y), \ y\in\Z^d_{0+},
\\
\infty, \ \quad\qquad y\notin\Z^d_{0+},
\end{cases}
\end{gather}
for all $x\in\Z^d_{0+}$.
 For  $t,k\ge0$ and $x\in\Z^d_{0+}$, let 
\begin{gather}                                                                                \label{D13}
 h_t^{(k)}(l,x):=\inf_{0\le j< \alpha_t(l)}[H^{(k)}(x+S_{t,t+j})-L_{t,t+j}(x+S_{t,t+j})].
\end{gather}

The function $H^{(k)}(y)$ 
describes the 
environment 
which is virgin for all $y\in\Z^d_{k+}$ and which 
 has no restrictions on the number of visits to  all states $ y\notin\Z^d_{k+}$.
The  function $h_t^{(k)}(l,x)$ 
 describes the behaviour in this environment of a random walk that starts at time $t\ge0$ from the state $x$.  Inequality (\ref{q0}) below shows that this environment has characteristics that
dominate the corresponding characteristics of any of our initial environments.

Note that 
\begin{gather}                                                                                \label{D+}
\{ h_t^{(k)}(l,x)= h_t(l,x),\,s_t(l)\ge0\}\subset \{\sup_{0\le j< t} S_j[1]<k\}.
\end{gather}
We use symbol $\infty$ in place of $0$ in \eqref{D12} because we like to use in Section 5 the following result (with ${\mathbf P} (\cdot)={\mathbf P}_+(\cdot)$):.

\begin{lemma}                                                                                                            \label{L1}
Under the assumptions $(A1)$ -- $(A3)$ and for each fixed $l\ge0$, given the event 
$\{S_t=x\}$ occurres, the joint conditional distribution of the random variables 
from the following family
\begin{gather*}                                                                                                    \label{D14}
\alpha_t(l),\ s_t(l),\ \vec S_{t,t+\alpha_t(l)},\  h^{(x[1])}_t(l,x),\,  ; \xi_{t+j},\,j\ge1
\end{gather*}
does not depend on $t\ge0$ and on $x\in\Z^d_{0+}$. In particular,
for all ${\cal C}\subset\Z^{d}_*$  
\begin{gather}                                                                                                    \label{D15}
{\mathbf P}(\vec S_{t,t+\alpha_t(l)}\in {\cal C},\, h^{(x[1])}_t(l,x)\ge0,\,s_t(l)\ge0|S_t=x)
\\                               \nonumber
={\mathbf P}(\vec S_{0,\alpha_0(l)}\in {\cal C},\,h_0(l,0)\ge0,\,s_0(l)\ge0)
={\mathbf P}_0(\alpha(l)<T_*,\,\vec S_{0,\alpha(l)}\in {\cal C}).
\end{gather}
\end{lemma} 
\proof 
The { first assertion} follows directly 
 from assumptions (A1)-(A3) and, in particular, from the time/space homogeneity of the random walk and from the homogeneity of the random environment in the positive half-space $\Z^d_{0+}$.
{To get \eqref{D15} we use \eqref{D+} too.}
\qed

\subsection{Auxiliary Lemmas}
Suppose that a 
random variable $T\ge0$ is such that
\begin{gather}                                                                                                       \label{w0} 
\{T\ge0\}=\cup_{t=0}^\infty\{T=t,\,S_T\in\mathbb X(t)\}
\quad \mbox{for\ some } \quad\mathbb X(t)\subset\Z^d.
\end{gather}
For a fixed $l>0$ and arbitrary sets ${\cal A},{\cal C}\subset\Z^d_*$, consider the event
\begin{gather}                                                                                                     \label{w1}
\tilde D:=\{T+\alpha_T(l)<T_*,\,\vec S_T\in {\cal A},\,
\vec S_{T,T+\alpha_T(l)}\in {\cal C}\}.
\end{gather}
Using (\ref{d2}) and (\ref{D5}), we may represent (\ref{w1}) in  the form
\begin{gather*}                                                                                                     \label{w2}
\tilde D=\{T<\infty,\, h(T-1)\ge0,\, h_T(l,S_T)\ge0,\,
\vec S_T\in {\cal A},\,
\vec S_{T,T+\alpha_T(l)}\in {\cal C}\}.
\end{gather*}
For  fixed  $t\ge0$ and $x\in\Z^d$, introduce events 
\begin{gather*}                                                                                                       \label{w3} 
\tilde A_{t,x}:=\{T=t,\,h(t-1)\ge0,\,\vec S_{t}\in {\cal A},\, S_t=x\},\ \ 
\tilde C_{t,x}:=\{h_{t}(l,x)\ge0,\, \vec S_{t,\alpha(n)}\in {\cal C}\}.
\end{gather*}
Clearly,
\begin{equation}                                                                                                   \label{w5}
{\mathbf P}_*(\tilde D)=
\sum_{t=0}^\infty \sum_{x\in \mathbb X(t)}{\mathbf P}_*(\tilde A_{t,x}\cdot\tilde  C_{t,x}). 
\end{equation}

Thus, we have the following elementary
\begin{lemma}  
Suppose that a random variable $T\ge0$ satisfies condition \eqref{w0}.                                                                                                         \label{L0}
Then 
 for  all ${\cal A},{\cal C}\subset\Z^d_*$ and each  $l>0$
equality \eqref{w5} takes place. 
In addition, 
if for all $t\ge0$ and $x\in\mathbb X(t)$ events $\tilde A_{t,x}$ and $\tilde C_{t,x}$ are pariwise independent and ${\mathbf P}_*(\tilde  C_{t,x})$ does not dependent 
on $t\ge0$ and $x\in\mathbb X(t)$, then we have
\begin{equation}                                                                                                   \label{w7}
{\mathbf P}_*(\tilde D)=
{\mathbf P}_*(T<T_*,\,\vec S_{T}\in {\cal A})
\cdot
{\mathbf P}_*(\alpha(l)<T_*,\, \vec S_{0,\alpha(l)}\in {\cal C}). 
\end{equation}
\end{lemma}
One can observe that the  sequence 
$\{S_t,\,H_{t}(x):x\in\Z^d\},\quad t=0,1,2,\dots,$ of infinite-dimensional random variables
forms an  
infinite-dimensional Markov chain.
In the proofs below we apply Lemma \ref{L0} four times for stopping times  $T\ge0$ of this Markov chain.

\begin{lemma}                                                                                                            \label{L2}
Under the assumptions $(A1)$ -- $(A3)$,
\begin{gather}                                                                                                      \label{d12}
{\mathbf P}_*(\alpha(n)<T_*,\,\vec S_{\alpha(k)}\in {\cal A},\,s(\alpha(k),\alpha(n))\ge0,\,\vec S_{\alpha(k),\alpha(n)}\in {\cal C})
\\                                                                                                        \nonumber
={\mathbf P}_*(\alpha(k)<T_*,\,\vec S_{\alpha(k)}\in {\cal A})\cdot
{\mathbf P}_0(\alpha(l)<T_*,\,
\vec S_{0,\alpha(l)}\in {\cal C})
\end{gather}
for  all ${\cal A},{\cal C}\subset\Z^d_*$ and each $n> k\ge0$ (where $l:=n-k>0$).
\end{lemma}

\proof We will apply Lemma \ref{L0} with $T=\alpha(l)$ and $\X(t)=\Z^d_{k}$. We have from 
 Lemma \ref{L1} that probability ${\mathbf P}_*(C_{t,x})$ does not depends on $t\ge0$ and $x\in \Z^d_{k}$. Hence, by~\eqref{D+}
\begin{gather}                                                                                                    \label{d18}
{\mathbf P}_*( \tilde C_{t,x})={\mathbf P}_*( \tilde C_{0,0})
={\mathbf P}_0(\alpha(l)<T_*,\,
\vec S_{0,\alpha(l)}\in {\cal C}).
\end{gather}

For any fixed $t\ge0$ and $x\in\Z^d_k$,  
random variables $\alpha(k)$, $\vec S_t$ and $h(t-1)$ 
that define the event $\tilde A_{t,x}$
are functions    only of the variables from the following two families:
\begin{gather}                                                                                                    \label{d11+-}
\{\xi_j: j\le t\}
\quad \mbox{and} \quad
\{H(y): y\notin\Z^d_{k+}\}.
\end{gather}
On the other hand, all random variables that determine the event $\tilde C_{t,x}$,  
are functions only of random variables from the  following two families:
\begin{gather}                                                                                                    \label{d11+}
\{\xi_j: j>t\}
\quad \mbox{and} \quad
\{H(y): y\in\Z^d_{k+}\}.
\end{gather}
Since the families in (\ref{d11+}) and (\ref{d11+-}) do not overlap, they are  independent.
Hence, events $\tilde A_{t,x}$ and $\tilde C_{t,x}$ are independent too.
This fact, together with \eqref{d18}, allows us to apply Lemma \ref{L0} to get \eqref{d12}.
\qed

\begin{lemma}                                                                                                            \label{L3}
Under the assumptions $(A1)$ -- $(A3)$,
\begin{gather}                                                                                                    \label{d23}
{\mathbf P}_*(\alpha(n)<T_*,\,\vec S_{\alpha(k)}\in {\cal A},\,\varkappa_*(n)=k,\,\vec S_{\alpha(k),\alpha(n)}\in {\cal C})
\\                                                                                                        \nonumber
={\mathbf P}_*(\alpha(k)<T_*,\,\vec S_{\alpha(k)}\in {\cal A})\cdot
{\mathbf P}_0(\varkappa_*(n-k)=0,\,\vec S_{0,\alpha_0(n-k)}\in {\cal C})
\end{gather}
for any $n> k\ge0$ and all ${\cal A},{\cal C}\subset\Z^d_*$.
\end{lemma}
\proof 
For each $n\ge1$ introduce the following subset of $\Z^d_*$:
\begin{gather}                                                                              \label{d21}
 {\cal C}_n^+:=\{(y_1,y_2,\dots)\in B_n^+:
\varkappa(n|\vec y)=0\quad \mbox{for}\quad \vec y=(0,y_1,y_2,\dots)\}.
\end{gather} 
We  assume in (\ref{d21}) that $y_0=0$ to avoid problems with the definition of the value $\alpha(n|\vec y)$.
It follows from (\ref{d21})  that
\begin{gather}                                                                                                      \label{d22}
\{\alpha(n)<T_*,\varkappa_*(n)=k\}
=\{\alpha(n)<T_*,\,s(\alpha(k),\alpha(n))\ge0,\,\vec S_{\alpha(k),\alpha(n)}\in {\cal C}^+_{n-k}\},
\end{gather}
where in (\ref{d22}) we used also that $\{\varkappa_*(n)=k\}\subset\{s(\alpha(k),\alpha(n))\ge0\}$.

If we compare now (\ref{d23}) and (\ref{d22}) with (\ref{d12}),
we can observe that (\ref{d23}) is a particular case of (\ref{d12}),
give that we replace in (\ref{d12})  ${\cal C}$ by ${\cal C}\cap{\cal C}^+_{n-k}$.
\qed

\subsection{Proof of Property 1 }
The first assertion of Property 1 
immediately follows from Lemma \ref{L3} with ${\cal A}={\cal C}=\Z^d_*$ since, in this case, we have from (\ref{d23}) that
\begin{gather}                                                                                                          \label{d24}
{\mathbf P}_*(B_n, \varkappa_*(n)=k)={\mathbf P}_*(\alpha(n)<T_*,\,\varkappa_*(n)=k)
\\                                                                                                        \nonumber
={\mathbf P}_*(\alpha(k)<T_*)\cdot {\mathbf P}_0(\varkappa_*(n-k)=0)
={\mathbf P}_*(B_k) \cdot{\mathbf P}_0 (\varkappa_*(n-k)=0).
\end{gather}
Since $B_n=\{\eta_*(n)\ge0\}$  and $\{\eta_*(n)\ge 1\}=\{\varkappa_*(n)\ge0\}$ by (\ref{a17}), we have, for $n=1,2,\ldots$,
\begin{gather}                                                                                                          \label{d24+}
{\mathbf P}_*(B_n)={\mathbf P}_*(\eta_*(n)=0)+{\mathbf P}_*(\eta_*(n)\ge1)
={\mathbf P}_*(\eta_*(n)=0)+\sum_{k=0}^{n-1} {\mathbf P}_*(\alpha(n)<T_*,\,\varkappa_*(n)=k).
\end{gather}
Thus, Property  1 follows from (\ref{d24}) and (\ref{d24+}). 



\section{Proof  of Theorem 1}
We will use functions $H^{(k)}(y)$ and $h_t^{(k)}(l,x)$ introduced in \eqref{D12} and \eqref{D13},
that have been already applied in Lemma \ref{L1}. 
These functions have the following useful properties:
\begin{gather}                                                                                \label{q0}
 H^{(k)}(y) \ge H_t(y) \quad \mbox{and} \quad
h_t^{(k)}(l,x)\ge  h_t(l,x) \quad   \forall y\in\Z^d,\ \forall x\in\Z^d_{0+}, \forall\ t,l\ge0 .
\end{gather}

\subsection{Main Lemma}
We are going to prove
\begin{pro}                                                                                                            \label{proA}
Under the assumptions $(A1)$ -- $(A4)$, there exists a constant $C<\infty$ such that
\begin{equation}                                                                                                   \label{Q1}
\forall\ n\ge0\quad 
{\mathbf P}_+(B_n)\le C {\mathbf P}_0(B_n).
\end{equation}
\end{pro}
The  proof 
is based on
several lemmas.
Introduce the following stopping time:
\begin{gather*}                                                                                                        \label{Q2}
\rho:=\inf\{t>0: S_t[1]=0\quad \mbox{but} \quad S_t\neq0\}\le\infty.
\end{gather*}
So  $\rho$ is the time of the first return to level $0$ by the
first component of our random walk, given that at least one of other coordinates differs from $0$.

\begin{lemma}                                                                                                            \label{M2}
For any $n>0$
\begin{equation}                                                                                                   \label{Q3}
P_\rho:={\mathbf P}_+(\rho<\alpha(n)< T_*)\le {\mathbf P}(\rho<\infty)\cdot {\mathbf P}_+(B_n).
\end{equation}
\end{lemma} 
\proof 
It follows from \eqref{D5} that
\begin{gather*} 
\{\rho<\alpha(n)< T_*\}=\{\rho<\infty,\,\alpha_\rho(n)<\infty, \,h(\rho-1)\ge0,\,
h_\rho(n,S_\rho)\ge0\}.
\end{gather*}
Since $h_\rho(n,S_\rho)\le h^{(0)}_\rho(n,S_\rho)$ by \eqref{q0}, we have

\begin{gather}                                                                                                       \label{Q4}
\{\rho<\alpha(n)< T_*\}\subset  
D:=\{\rho<\infty,\,\alpha_\rho(n)<\infty, \, h^{(0)}_\rho(n,S_\rho)\ge0\}.
\end{gather}
Introduce the events 
\begin{gather*}                                                                                                       \label{Q5}
A_{t,x}:=\{\rho=t,\,S_t=x\},\quad
C_{t,x}:=\{\alpha_t(n)<\infty,\, h^{(0)}_t(n,x)\ge0\}.
\end{gather*}
By Lemma \ref{L1},  probability ${\mathbf P}(C_{t,x})$ does not depends on $t\ge0$ and $x\in \Z^d_{0}$. Hence, by~\eqref{D+} 
\begin{gather}                                                                                                      \label{Q7}
{\mathbf P}_*(C_{t,x})={\mathbf P}_*(\alpha_0(n)<\infty,\,  h^{(0)}_0(\alpha_0(n),0)\ge0)
={\mathbf P}_+(\alpha_0(n)<T_*)={\mathbf P}_+(B_n),
\end{gather}
since $\alpha_0(n)=\alpha(n)$.

Now we  apply Lemma \ref{L0}  with  
$T=\rho$ and $\X(t)=\X_0:=\Z^d_{0}\setminus\{0\}$,
and with $ h^{(k)}_k(l,x)$ in place of $h_k(l,S_k)$.
For fixed values $t>0$ and $x\in \Z^d_{0}$,
random variables $\alpha_t(n)$ and $h_t(n,x)$
are functions only of random variables  from \eqref{d11+} with $k=0$,
since $H(y)=\infty $ for all $y\notin \Z^d_{k+}$.

On the other hand, event $A_{t,x}$ 
does not depend on the environment and is determined by   
the variables $\{\xi_j: j\le t\}$.
Hence, events $\tilde A_{t,x}$ and $\tilde C_{t,x}$ do not  depend on each other, and we may apply Lemma \ref{L0}. Using also \eqref{Q7} and \eqref{Q4}, we obtain 
\begin{equation*}  
{\mathbf P}_*(\rho<\alpha(n)< T_*)\le {\mathbf P}_*(D)=
\sum_{t=1}^\infty \sum_{x\in \Z^d_{0}\setminus\{0\}}{\mathbf P}(A_{t,x})
=\sum_{t=1}^\infty{\mathbf P}(\rho=t, S_t\neq0)={\mathbf P}(\rho<\infty).
\end{equation*}
Thus \eqref{Q3} is proved.
\qed

{According to \eqref{a9} } introduce the following stopping times:
\begin{gather*}                                                                                                        \label{Q13}
\rho_0=0\quad \mbox{and} \quad
\rho_i:=\inf\{t>\rho_{i-1}: S_t=0\}\le\infty,
\quad i=1,2,\dots .
\end{gather*}
So $\rho_i$ is the time of the $i$-th return to $0$ of our random walk.
It is easy to see that, for any  $n>0$,
\begin{equation}                                                                                                   \label{Q14}
{\mathbf P}_+(B_n)\le{\mathbf P}_+(\rho<\alpha(n)< T_*)+\sum_{i=0}^\infty {\mathbf P}(D_i),
\end{equation}
where
\begin{equation*}                                                                                                   \label{Q15}
D_i=D_i(n):= \{\rho_i<\alpha(n)<\min (\rho_{i+1}, \rho )\le\infty,\,\alpha(n)< T_*\}. 
\end{equation*}

\begin{lemma}                                                                                                            \label{M3}
For any $n>0$
\begin{gather}                                                                                                        \label{Q16}
{\mathbf P}_+(D_i)\le {\mathbf P}(\rho_i<\infty)
\cdot {\mathbf P}( H(0)>i)\cdot {\mathbf P}_0(B_n). 
\end{gather}
\end{lemma} 
\proof
Underline that, on the event $\{\rho_i<\alpha(n)<\min (\rho_{i+1},\rho)\}$, we have $s(\rho_i,\alpha(n))>0$, due to the skip-free property of the random walk.  Thus
\begin{equation}                                                                                                   \label{Q17}
D_i\subset \hat D_i:=
 \{\rho_i<\alpha(n)
< T_*,\,s(\rho_i,\alpha(n))>0\}. 
\end{equation}
Since $S_{\rho_i}=0$, we have from  \eqref{D5} that
\begin{gather}                                                                                                       \label{Q18}
\hat D_i=\{\rho_i<\infty,\,\rho_i<\rho,\,\alpha_{\rho_i}(n)<\infty, \,h(\rho_i-1)\ge0,\,
h_{\rho_i}(n,0)\ge0\}.
\end{gather}
Since $H_{t+j}(y)\le H_{j}(y)$ for all $y\in\Z^d$ and $t,j\ge0$, we have from \eqref{d3} that 
\begin{gather}                                                                                                \label{Q19}
h_t(n,0)=\inf_{0\le j< \alpha_t(n)}H_{t+j}(S_{t,t+j})
\le \tilde h_t(n,0):=\inf_{0\le j< \alpha_t(n)}H_{j}(S_{t,t+j})
\end{gather}
for all possible $t\ge0$. Note also that $h(\rho_i)\le H_{\rho_i}(0)$.
This fact and (\ref{Q17}), (\ref{Q18}) and (\ref{Q19}) with $t=\rho_i$ yield
\begin{gather*}                                                                                                       \label{Q20}
\hat D_i\subset \tilde D_i:=
\{\rho_i<\infty,\,\rho_i<\rho,\,\alpha_{\rho_i}(n)<\infty, \,H_{\rho_i}(0)\ge0,\,
\tilde h_{\rho_i}(n,0)\ge0,\,s(\rho_i,\alpha(n))>0\}.
\end{gather*}

Introduce the events 
\begin{gather}                                                                                                       \label{Q21}
A_{i,t}:=\{\rho>\rho_i=t,\,H_t(0)\ge0\},\quad
C_{t}:=\{\alpha_t(n)<\infty,\,\tilde h_t(n,0)\ge0,\,s(t,\alpha(n))>0\}.
\end{gather}
Comparing definition \eqref{Q19} with that in \eqref{d3} and \eqref{D3}, we can see 
that  $\tilde h_t(n,0)= h_t(n,0)$ and that probability ${\mathbf P}_*(C_{t})$ does not depend on $t\ge0$. Hence, 
\begin{gather}                                                                                                      \label{Q23}
{\mathbf P}_*(C_{t})={\mathbf P}_*(C_{0})={\mathbf P}_*(\alpha_0(n)<\infty,\,  \tilde h_0(n,0)\ge0,\,s(0,\alpha(n))>0)
\\                                    \nonumber
\le {\mathbf P}_*(\alpha_0(n)<\infty,\,  h_0(n,0)\ge0,\,s(0,\alpha(n))\ge0)
={\mathbf P}_0(\alpha_0(n)<T_*)={\mathbf P}_0(B_n),
\end{gather}
because $\alpha_0(n)=\alpha(n)$. We have also from \eqref{Q21} that
\begin{gather}                                                                                                        \label{Q25}
{\mathbf P}_*(A_{i,t})
={\mathbf P}( \rho>\rho_i=t){\mathbf P}(H_t(0)\ge0)
={\mathbf P}( \rho>\rho_i=t){\mathbf P}( H(0)>i),
\end{gather}
because $H_{\rho_i}(0)=H_0(0)-i=H(0)-i-1$. 

Now we apply Lemma \ref{L0} with $T=\rho_i$ and $\X(t)=\{0\}$,
and with $ \tilde h_t(n,0)$ in place of $ h_t(n,0)$.
Now note that, for each value $t>0$, under the condition $s(t,\alpha(n))>0$, the 
 random variables $\alpha_t(n)$, $s(t,\alpha(n))$ and $h_t(n,x)$ 
(which determine event $C_t$)
are functions only of random variables from \eqref{d11+}, with $k=1$.
On the other hand, event $A_{i,t}$  is defined by    
the variable $H_t(0)$ and by the family $\{\xi_j: j\le t\}$.
Hence, the events $A_{i,t}$ and $C_{t}$ are  independent  
and we may apply Lemma \ref{L0} again. Using also  \eqref{Q23} and \eqref{Q25}, we get: 
\begin{equation*}  
{\mathbf P}_+(D_i)\le{\mathbf P}_+(\hat D_i)\le{\mathbf P}_+(D_i)
=\sum_{t=1}^\infty {\mathbf P}(A_{i,t}){\mathbf P}(C_{t})
\le\sum_{t=1}^\infty {\mathbf P}(\rho>\rho_i=t){\mathbf P}( H(0)>i){\mathbf P}_0(B_n)
.
\end{equation*}
So, inequality (\ref{Q16}) follows.
\qed

Introduce the notation
\begin{gather*}                                                                                                   \label{Q27}
p^*:={\mathbf P}(\rho<\infty),\quad
p_1:={\mathbf P}(\rho_1<\infty,\,\rho_1<\rho),
\\                                                                                                  \label{Q28}
C^*:=\sum_{i=0}^\infty p_1^i{\mathbf P}( H(0)>i)\le1+{\mathbf E} H(0).
\end{gather*}
Substituting the results of Lemmas \ref{M2} and \ref{M3}   into (\ref{Q14}), we obtain
$$
1-{\mathbf P}_+(B_n)\le C^*{\mathbf P}_0(B_n)+p^*{\mathbf P}_+(B_n).
$$
Thus, under the assumptions $(A1)$ -- $(A3)$ ,
\begin{equation}                                                                                                   \label{Q33}
(1-p^*){\mathbf P}_+(B_n)\le C^*{\mathbf P}_0(B_n).
\end{equation}

One can easily conclude that, under any of assumptions (a)-(c)  
in $(A4)$,
the following inequalities hold: 
\begin{equation}                                                                                                     \label{Q35}
p_*<1  \qquad\mbox{and} \qquad
C_*<\infty.
\end{equation}
Here is the only place in the paper where 
the assumption $(A4)$ is used.

From (\ref{Q33}) and (\ref{Q35}) we obtain the assertion of Property \ref{proA}]
with $C=C^*/(1-p^*)$.

Note that for $n=0$ inequality \eqref{Q1} follows from \eqref{a5} since $C\ge1$.
\qed

\subsection{Using Submultiplicativity }
  In this  Subsection we prove first that                                                                                                       
  \begin{gather}                                                                                                        \label{q51}
\forall k,l\ge0\quad{\mathbf P}_*(B_k) {\mathbf P}_0(B_l) 
\le{\mathbf P}(B_{k+l}) \le {\mathbf P}(B_k) {\mathbf P}_+(B_l) .
\end{gather}
 Using this 
form of sub/supermultiplicativity we show that
\begin{gather}                                                                                                            \label{q52}
1\ge q_+:=\inf_{n\ge1} \sqrt[n]{{\mathbf P} _+(B_n)}
=  q:=
\sup_{n\ge1} \sqrt[n]{{\mathbf P} _-(B_n)}
\ge {\mathbf P} (\xi_1[1]=1)>0.
 \end{gather}
After that,  we prove the following
\begin{pro}                                                                                                          \label{proB}
Under the assumptions $(A1)$ -- $(A4)$, 
relations \eqref{q51} and  \eqref{q52}  take place. Moreover
\begin{gather}                                                                                 \label{q53}
0< {\mathbf P}_*(B_0)/C\le {\mathbf P}_* (B_n)/q^{n}
\le C{\mathbf P}_*(B_0) \le C<\infty
\quad\forall\ n\ge0 .
\end{gather}
 \end{pro} 
  Note that for $l=0$ inequality \eqref{q51}  
immediately follows from \eqref{a5}. 
We prove now the following lemma.

\begin{lemma}                                                                                                            \label{M6}
Under the assumptions $(A1)$ -- $(A3)$, inequality \eqref{q51} takes place for all $k\ge0$ and $l>0$.
\end{lemma} 
\proof 
Applying Lemma \ref{L2} with $n=k+l$, we get 
\begin{gather*}                                                                                                      \nonumber
{\mathbf P}_*(B_{k+l})={\mathbf P}_* (\alpha(k+l)< T_*)
\ge{\mathbf P}_* (\alpha(l+k)< T_*, \,  s(\alpha(k),\alpha(k+l))\ge0)
\\                                                                                                                            \label{q55}
={\mathbf P}_* (\alpha(k)< T_*)
\cdot{\mathbf P}_0 (\alpha(l)< T_*) ={\mathbf P}_*(B_k) {\mathbf P}_0(B_l).
\end{gather*}
This is the first inequality in \eqref{q51}.

Next, it follows from \eqref{D5} that
\begin{gather*} 
B_{k+l}=\{\alpha(k+l)< T_*\}
=\{\alpha(k)<\infty,\,\alpha_k(l)<\infty, \,h(k-1)\ge0,\,h_k(l,S_k)\ge0\}.
\end{gather*}
Since $h_k(l,S_k)\le h^{(k)}_k(l,S_k)$ by \eqref{q0}, we have
\begin{gather*}                                                                                                       \label{Q54}
B_{k+l}\subset  \tilde D:=\{\alpha(k)<\infty,\,\alpha_k(l)<\infty, \,h(k-1)\ge0,\, h^{(k)}_k(l,S_k)\ge0\}.
\end{gather*}

Now we apply Lemma \ref{L0} with the same 
$T=\alpha(l)$ and $\X(t)=\Z^d_{k}$ as in the proof of Lemma \ref{L2},
but with $ h^{(k)}_k(l,x)$ in place of $h_k(l,S_k)$. 
Introduce events 
\begin{gather*}                                                                                                       \label{Q55}
\tilde A_{t,x}:=\{\alpha(k)=t,, \,h(t-1)\ge0,\,S_t=x\},\quad
\tilde C_{t,x}:=\{\alpha_k(l)<\infty,\, h^{(k)}_k(l,x)\ge0\}.
\end{gather*}
By Lemma \ref{L1},  probability ${\mathbf P}(\tilde C_{t,x})$ does not depends on $t\ge0$ and $x\in \Z^d_{k}$,
\begin{gather}                                                                                                      \label{Q57}
{\mathbf P}_*(\tilde C_{t,x})={\mathbf P}_*(\tilde C_{0,0})={\mathbf P}_*(\alpha_0(l)<\infty,\,
 h^{(0)}_0(\alpha_0(l))\ge0)
={\mathbf P}_+(\alpha_0(l)<T_*)={\mathbf P}_+(B_l),
\end{gather}
since $\alpha_0(l)=\alpha(l)$  for $l>0$. 

Now, for fixed values $t>0$ and $x\in \Z^d_{k}$,
 random variables $\alpha_t(n)$ and $h_t(n,x)$,
which define event $\tilde C_{t,x}$,
are functions only of random variables from \eqref{d11+},
since $H(y)=\infty $ for all $y\notin \Z^d_{k+}$.
On the other hand, event $\tilde A_{t,x}$ 
is defined by  random variables $\alpha(k)$, $h(t-1)$ and~$S_t$
which are functions  of  the variables from \eqref{d11+}.
Hence, events $\tilde A_{t,x}$ and $\tilde C_{t,x}$ are independent and we can apply Lemma \ref{L0}. Using also \eqref{Q57}, we obtain 
${\mathbf P}_*(B_{k+l}) \le {\mathbf P}_*(\tilde D)= {\mathbf P}_*(B_k) {\mathbf P}_+(B_l)$
as a  result.

Thus,  second inequality in \eqref{q51} is proved.
\qed

\proof[Proof of Property \ref{proB}]
Using probabilities ${\mathbf P}_0(\cdot)$ and ${\mathbf P}_+(\cdot)$
instead of ${\mathbf P}_*(\cdot)$,
we have from \eqref{q51} and \eqref{a5}
 that, for all $k,l\ge1$,  
\begin{gather*} 
{\mathbf P}_0(B_k) {\mathbf P}_0(B_{kl-k}) 
\le{\mathbf P}_0(B_{kl})\le{\mathbf P}_+(B_{kl}) \le {\mathbf P}_+(B_{kl-l}) {\mathbf P}_+(B_l) .
\end{gather*}
Then the induction argument leads to 
\begin{gather}                                                                                                        \label{q57}
({\mathbf P}_0)^l(B_k)  
\le{\mathbf P}_0(B_{kl})\le{\mathbf P}_+(B_{kl}) \le ( {\mathbf P}_+)^k(B_l) .
\end{gather}
Taking the $k$th root of the both sides of inequality \eqref{q57}, 
we arrive to
\begin{gather}                                                                                                            \label{q58}
 \forall\ k,l\ge1\quad 
\sqrt[k]{{\mathbf P} _0(B_k)}
\le \sqrt[l]{{\mathbf P} _+(B_l)}.
 \end{gather}
Taking in \eqref{q58} supremum in  $k\ge1$ and infimum in $l\ge1$,
we obtain $q\le q_+$.

On another hand, from \eqref{Q1} and the definition of $q_+$ in \eqref{q52}, we have
\begin{equation}                                                                                                               \label{q59}
q_+^n\le{\mathbf P}_+(B_n)
\le C {\mathbf P}_0(B_n)\le  Cq^n.
\end{equation}
Hence, $q_+\le \sqrt[n]{C} q\to q$. So, we proved that $q_+\le q$ and hence \eqref{q52} follows 
from  \eqref{a3} with ${\mathbf P} _0(B_n)\ge{\mathbf P}^n(\xi_1[1]=1)$. 

Next, it follows from \eqref{q51} and \eqref{Q1}   
with $k=0$ and $l=n$ that 
\begin{gather}                                                                                                        \label{q62}
{\mathbf P}_*(B_{n}) \le {\mathbf P}_*(B_0) {\mathbf P}_+(B_n)
\le C{\mathbf P}_*(B_0) {\mathbf P}_0(B_n)\le  C{\mathbf P}_*(B_0) q^n\le Cq^n.
\end{gather}
Here we also used \eqref{q59}. On the other hand, using again \eqref{q51},  \eqref{Q1} and \eqref{q59},
we get
\begin{gather}                                                                                                        \label{q63}
{\mathbf P}_*(B_n) \ge {\mathbf P}_*(B_0) {\mathbf P}_0(B_{n}) \ge 
{\mathbf P}_*(B_0) {\mathbf P}_+(B_n) /C\ge {\mathbf P}_*(B_0) q_+^n/C .
\end{gather}
Now, all inequalities in \eqref{q53} follow from \eqref{q62} and \eqref{q63}.
\qed

 
\subsection{Proof of Theorem 1}  
With $q$ from (\ref{q52}), introduce the following notation:
\begin{gather}                                                                                                          \label{b71}
a_n:=\frac{{\mathbf P}_{0}(\varkappa_*(n)=0)}{q^n},\ 
 b_n:=\frac{{\mathbf P}_*(\eta_*(n)=0)}{q^n},\ 
u_n:=\frac{{\mathbf P}_*(B_n)}{q^n},\  v_n:=\frac{{\mathbf P}_0(B_n)}{q^n}.
\end{gather}
Multiplying equalities \eqref{a20} and \eqref{a21} by $q^{-n}$, we obtain for all $n\ge1$ that
\begin{gather}                                                                                                          \label{b73}
u_n=b_n+\sum_{k=0}^{n-1} u_ka_{n-k}=b_n+\sum_{l=1}^{n}a_l u_{n-l},
\\                                                                                                         \label{b74}
v_n=\sum_{l=1}^{n} a_lv_{n-l},  \quad \mbox{where} \quad
v_0=1
\quad \mbox{and} \quad a_1>0.
\end{gather}
The last property in \eqref{b74} follows from \eqref{b5}.

{We have  from \eqref{b71},} \eqref{q52} and \eqref{q53} that
\begin{gather}                                                                                                          \label{b72}
u_0=b_0={\mathbf P}_*(B_0)>0, \ \ 0<u_0/C\le u_n\le C<\infty,\ \ 0<1/C\le v_n\le1
\ \ \forall n\ge1 .
\end{gather}
In addition, we have  from \eqref{a15}, \eqref{b5} and \eqref{b72} that
\begin{gather}                                                                                                          \label{b77}
v_1=a_1>0,\quad
0\le a_n\le v_n
\quad \mbox{and} \ \
0<v_1^n\le v_n\le1 \ \ \forall n\ge1 .
\end{gather}

There are two possible scenarios, either $a_n<1$ for all $n$ or
$a_M=1$ for some $M\ge 1$. We start with the latter case which is, in fact, degenerative. 

\begin{lemma}                                                                                                            \label{M9}
If $a_{M}=1$ for some $M\ge1$, 
then $M=1$ and 
the assertions of Theorem 1 do hold with $q={\mathbf P}_0(\varkappa_*(1)=0)$.
\end{lemma} 
\proof Since $1=a_M\le v_M\le1$ by (\ref{b77}), we have  $v_M=1$.
Then, by (\ref{b74}),  
$$
v_M-a_M=0=\sum_{l=1}^{M-1} a_lv_{n-l}\ge a_1v_{M-1}>0
\quad \mbox{if} \quad M>1.
$$
So we must have $M=1$. Then $v_1=a_1=1=v_1^n\le v_n\le 1$ for all $n\ge1$ by (\ref{b77}). 
Hence, $v_n=1$ for all $n\ge1$ and, by (\ref{b74}), 
$$
v_n-a_1v_{n-1}=1-1=0=\sum_{l=2}^{n} a_lv_{n-l}=\sum_{l=2}^{n} a_l
\quad \mbox{when} \quad n\ge2.
$$
Thus, $a_l=0$ for all $l\ge2$ and equation (\ref{q1}) reduces to 
${\mathbf P}_0(\varkappa_*(1)=0)/q=1$.
Hence all assertions of Theorem 1 hold with $q={\mathbf P}_0(\varkappa_*(1)=0)$.
\qed

Consider now the main case where
$0\le a_k<1$ for all $k\ge1$. It is known 
(see, for example, Section 13.4 in  the 1st Vollume of the Feller's book \cite{Fe})
that there are only four possibilities for the solutions to equation \eqref{b74}:
\\
(a)\quad $0<\alpha:=\sum_{k\ge1}a_k<1$ \ and \ $v_n\to0$;
\\
(b)\quad $\alpha=1$,\ $\mu=\sum_{k\ge1}ka_k=\infty$\ and \ $v_n\to0$;
\\
(c)\quad $\alpha=1$,\ $1\le\mu<\infty$\ and \ $v_n\to1/\mu>0$ since $a_1>0$;
\\
(d)\quad $\alpha\in(1,\infty]$\ and \ $v_n\to\infty$.
\\
It is easy to see that (c) is the  only  possibility which does not contradict to  inequalities \eqref{b72}.   
Hence, $\alpha=1$,\ $\mu<\infty$, and \eqref{q1} with \eqref{q2} follow.

Now, 
we again use \cite{Fe} to evaluate $\psi_0=\sum_{k\ge1}b_k$.
From \eqref{b72} and \eqref{b73} with $v_n\to1/\mu>0$ we obtain
\begin{gather*}                                                                                                         
C\ge u_n=\sum_{k=0}^{n} b_ku_{n-k}\to \sum_{k\ge0}b_k/\mu=\psi_0/\mu
\ge b_0/\mu={\mathbf P}_*(B_0)/\mu>0
\end{gather*}
by assumption ${(A2)}$.
 So, we obtain inequality \eqref{q3}  with \ $C\mu\ge\psi_0\ge{\mathbf P}_*(B_0)>0$.

Thus, Theorem 1 is proved.


\section{Proof of Theorem 2} 
 
We suppose that assumptions $(A1)$ --- $(A4)$ continue to hold.

We start with a few preliminary comments. 
If follows directly from \eqref{b3} and \eqref{b6} that, for any integers $K\ge k\ge 0$ and all vectors $\vec y_K=(y_0,\dots,y_K)\in\Z^{(K+1)\times d}$,
 \begin{gather}                                                                                                              \label{z1}
{\mathbf P}_* (\alpha(k)=K<T_*,\eta_*(k)=0,\vec S_K=\vec y_K)
=\psi_0q^k{\mathbf P} (\oo\nu_0 = k,\, \oo T_0=K,\,  \ww S_K=y_K)).
\end{gather}
Similarly, it follows from \eqref{b4} and \eqref{b7} that, 
for any integers $L\ge l\ge1$ and all $\vec x_L=(x_1,\dots,x_L)\in\Z^{L\times d}$
\begin{gather}                                                                                                              \label{z2}
  {\mathbf P}_{0} ( \alpha(l)=L<T_*,\varkappa_*(l)=0,\vec S_{0,L}=\vec x_L)
=q^l{\mathbf P} (\oo{\lambda}_0 = l, 
\oo \tau_0=L, \ww S_{0,L}=\vec x_L) 
\\                                                                                                      \nonumber
=q^l{\mathbf P} (\oo{\lambda}_m\equiv\oo{\nu}_m-\oo{\nu}_{m-1} = l, 
\oo \tau_m\equiv\oo{T}_m-\oo{T}_{m-1}=L, \ww S_{\oo{T}_{m-1},\oo{T}_{m}}=\vec x_L).
\end{gather} 

In the proof of the following lemma we 
repeat  in more detail the description of the core random sequence, introduced in subsection 3.3.
\begin{lemma}                                                                                                        \label{L6}
 Suppose that numbers $N\ge n\ge m\ge 0$ and vector $\vec y_N=(y_0,\dots,y_N)\in\Z^{(N+1)\times d}$ 
are such that
\begin{gather}                                                                                                \label{z3}
 \alpha(n|\vec y_N)=N\ge0 \quad \mbox{and} \quad \eta(n|\vec y_N)=m\ge0.
\end{gather}
Then
\begin{gather}                                                                                                               \label{z4}
{\mathbf P}_* (\alpha(n)=N<T_*,\eta(n)=m,\vec S_N=\vec y_N)
=\psi_0q^n{\mathbf P} (\oo\nu(m)=n,\oo{\alpha}(n)=N,\ww S_N=\vec y_N).
\end{gather}
Moreover, all random variables in (\ref{z4}) are  deterministic functions only of random variables from  the initial block
and from the first $m$ blocks in (\ref{b14}).
\end{lemma}
We will prove the lemma
by induction in $m$. For $m=0$,  \eqref{z4} follows from \eqref{z1} (with $k$ in place of $n$ and $K$ in place of $N$) that has been verified already. 

 Let $m$ be a strictly positive number and suppose that \eqref{z4} holds  for all possible $N$ and ${\vec y}_N$ in the case 
 $\eta_*(n)=m-1\ge0$. Now take the numbers and a vector satisfying \eqref{z3}.
  Then, for some integers $k$ and~$K$, 
\begin{gather}                                                                                                              \label{z5}
\varkappa(n|\vec y_N)=k\in[0,n-1]\quad \mbox{and} \quad\alpha(k|\vec y_N)=K\in[0,N-1].
\end{gather}
Let 
 \begin{gather}                                                                                                              \label{z6}
\vec y_K=(y_0,\dots,y_K),\quad \vec y_{K,N}=(y_{K+1}-y_K,\dots,y_N-y_K),\quad N>K\ge0.
\end{gather}
We have from (\ref{z5}) that 
\begin{gather}                                                                                                        \label{z7}
\{\alpha(n)<T_*,\eta(n)=m,\varkappa(n)=k,\}=\{\alpha(n)<T_*,\eta(k)=m-1,\varkappa(n)=k\}.
\end{gather}
Hence, by (\ref{z6}) and (\ref{z7}),
\begin{gather}                                                                                                         \label{z8}
{\mathbf P}_* (\alpha(n)=N<T_*,\eta(n)=m,\vec S_N=\vec y_N)
\\                                                                                                            \nonumber
=     
{\mathbf P}_* (\alpha(n)=N<T_*,\eta(n)=m,\varkappa(n)=k,\alpha(k)=K,\vec S_K=\vec y_K,
\vec S_{K,N}=\vec y_{K,N}) . 
 \end{gather}
 Now  we apply Lemma \ref{L3} with special sets 
${\cal A}=\{\vec y_{K}\}$ and ${\cal C}=\{\vec y_{K,N}\}$ containing only one trajectory each. Then
\begin{gather}                                                                                                         \label{z9}
{\mathbf P}_* (\alpha(n)=N<T_*,\eta(n)=m,\vec S_N=\vec y_N)
\\                                                                                                            \nonumber 
={\mathbf P}_* (\alpha(k)=K<T_*,\eta(k)=m-1,\vec S_K=\vec y_K)  
\cdot {\mathbf P}_0 (\alpha(l)=L<T_*,\varkappa(l)=0,\vec S_{0,L}=\vec y_{K,N}).
\end{gather}

 Clearly, $\eta_*(k)=m-1$ by (\ref{z7}). Hence, by the induction base, we have 
 that
\begin{gather}                                                                                                               \label{z10}
{\mathbf P}_* (\alpha(k)=K<T_*,\eta_*(k)=m-1,\vec S_K=\vec y_K)
=\psi_0q^k{\mathbf P} (\ww S_K=\vec y_K,\oo{\alpha}(k)=K,\oo\nu_{m-1}=k).  
\end{gather}
Now use (\ref{z2}) with 
\begin{gather}                                                                                                              \label{z11}
\oo{\nu}_m=n,\quad\oo{\nu}_{m-1}=\varkappa_*(\oo{\nu}_m)=k,\quad
\oo{T}_m=\oo\alpha(\oo\nu_m)=N,\quad\oo{T}_{m-1}=\oo\alpha(\oo\nu_{m-1})=K.
\end{gather} 
Let $l=n-k$, $L=N_K$ and $\vec x_L=\vec y_{K,N}$.
Substituting (\ref{z10}) and (\ref{z2}) into  (\ref{z9}), we obtain from (\ref{z9}) and  (\ref{z11})  that  
\begin{gather*}                                                                                            \nonumber 
{\mathbf P}_* (\alpha(n)=N<T_*,\eta_*(n)=m,\vec S_N=\vec y_N)
\\                                                                                                            \label{z12}
=\psi_0q^k{\mathbf P} (\ww S_K=\vec y_K,\oo{\alpha}(k)=K,\oo\nu_{m-1}=k) 
\\                                                                                                      \nonumber
\cdot q^l{\mathbf P} (\oo{\nu}_m-\oo{\nu}_{m-1} = l, 
\oo{T}_m-\oo{T}_{m-1}=L, \ww S_{\oo{T}_{m-1},\oo{T}_{m}}=\vec x_L).
\end{gather*}
Notice that the $m$-th block in (\ref{b14}) is independent of the previous ones. Hence,  (\ref{z12}) may be represented as 
\begin{gather*}                                                                                                 \nonumber
{\mathbf P}_* (\alpha(n)=N<T_*,\eta_*(n)=m,\vec S_N=\vec y_N)
\\                                                                                                            \label{z25}
=\psi_0q^{k+l}{\mathbf P} (\oo{\nu}_m=n,\oo{T}_m=\oo{\alpha}(\oo{\nu}_m)=N,\ww S_{\oo{T}_{m-1}}=\vec y_K,
 \ww S_{\oo{T}_{m-1},\oo{T}_{m}}=\vec x_L=\vec y_{K,N}) 
\\                                                                                                      \nonumber
=\psi_0q^n{\mathbf P} (\oo\nu(m)=n,\oo{\alpha}(n)=N,\ww S_N=\vec y_N)
\end{gather*}
Therefore, we have completed the induction step. This ends the proof of
Lemma \ref{L6}.
\qed

To prove Theorem 2, note that any
set ${\cal A}\in \Z^d_*$ may be represented as 
\begin{gather*}                                                                                                               \label{z27}
{\cal A}=\cup_{N=0}^\infty A_N, \quad \mbox{where} \quad A_N\subset \Z^{(N+1)\times d},\quad N=0,1,2,\dots.
\end{gather*}
So, all vectors $\vec y_N=(y_0,\dots,y_N)$ from $A_N$ satisfy (\ref{z3}).
Then summing up the LHS's and RHS's of
\eqref{z4} over $N$ and $\vec y_N\in  A_N$  leads to \eqref{b21}.

Thus, we have finished with the proofs of all our results.

\section{Remarks}

{\bf Remark 7.1.}
In our Assumptions (A1)-(A3), we assume that the environment in ``virgin'' only in a half-space and that the random walk starts either from the other half-space or from a boundary point. Here is a scenario that may lead to such situation.

Assume that at some time instant $-\infty \le-N<0$ in the past the whole environment in $Z^d$ was ``virgin'', i.e.all  the 
random variables $\{H_{-N}(x),\,x\in\Z^d\}$ were i.i.d. 
Assume that our random walk had started at time $t>-N$.
This assumption implies that
 \begin{gather}                                                                                                      \label{z16}
H(x)=H_{-1}(x)\le H_{-N}(x),\quad\forall\,x\in\Z^d.
\end{gather}
 We then assume that the  trajectory of our random walk  
on the time interval $-N<t<0$ is unobservable (it is the ``dark history''), and that we start to observe the trajectory only at time $t=0$ when we realize that the environment 
is still virgin in the half-space 
$\Z^d_{0+}$ (see \eqref{a0} for definition), so that
 \begin{gather*}  
H(x)=H_{-1}(x)= H_{-N}(x),\quad\forall\,x\in\Z^d_{0+}.
\end{gather*}
Thus we arrive to our model with $S_{-1}\notin\Z^d_{0+}$ (and, hence, with $S_0[1]\le0$). 
\\
Note that our condition~(\ref{a2}) is more general than 
(\ref{z16}). 

{\bf Remark 7.2.}
Here is a link to random walks conditioned not to leave a certain subspace. 
We may consider the trajectory 
 $S_0,S_1,\dots, S_{\alpha(n)}$ conditioned on the event that the first coordinate stays positive by time $\alpha(n)$, i.e.  
$B_n=\{\min_{0\le t<\alpha(n)} S_t[1]\ge0\}$.
Then, in our notation, the event $B_n$ may be represented as $B_n=\{\alpha(n) < T_*\}$ if we consider that ``extreme'' environment of the form: for $x=(x[1],\ldots,x[d])$, 
\begin{equation*}                                                                                                 \label{i16}
H(x)=0\ \ \mbox{when}\ \ x[1]<0,
\quad \mbox{and}\quad
H(x)=\infty\ \ \mbox{when}\ \ x[1]\ge0.
\end{equation*}
Thus, there is no restrictions on the upper half-space with $x[1]\ge0$,
and it is prohibited  to visit the  lower half-space with $x[1]<0$. 

Note that the case  ${\mathbf E}\xi_1[1]>0$  is simple, since here the initial sequence itself has a regenerative structure and
(\ref{i13}), (\ref{i12}) and (\ref{i10}) take place with $q=1$.
In the case ${\mathbf E} \xi_1[1]<0$, there is only one $q\in (0,1)$ that solves the equation
\begin{equation}                                                                                                   \label{z18}
\sum_{k=-1}^{\infty}q^k {\mathbf P} (\xi_1[1]=-k) = 1.
\end{equation}
Applying the corresponding exponential change of measure (the Cram\'er 
transform) to the distribution of $\xi_1$, we obtain 
(\ref{i13}), (\ref{i10}) and \eqref{i12} with $q<1$ from (\ref{z18}).

{\bf Remark 7.3}.
We may present a more detailed version of Theorem \ref{T-repr},
containing a formula that relates joint distributions of blocks from (\ref{a27}) and (\ref{a28})
 with independent blocks of the core process. 
Consider arbitrary numbers such that 
\begin{gather}                                                                                     \label{z19}
0\le L_0<\dots<L_m=n,\quad 0\le K_i<\dots<K_m, \quad\vec y\in\Z^{(K_0+1)\times d},
\\ \nonumber
1\le l_i:=L_i-L_{i-1}\le k_i:=K_i-K_{i-1},\quad \vec x_i\in\Z^{(K_i-K_{i-1})\times d},
 \quad\forall i=1,\dots, m.
\end{gather}
Below we use the notation for vectors introduced in  (\ref{a25}) and (\ref{b13}).
\begin{corollary}                                                                                                          \label{C4}
For any $ n=L_m\ge m\ge1$   and any numbers from (\ref{z19})
\begin{gather*}   \nonumber %
 {\mathbf P}_* ( \{\eta_*(n)=m, \nu_0(n)=L_0,\alpha(\nu_0(n))=K_0,\vec S_{K_0}=\vec y\}
\\                                                                                                                         \label{z20}
 \cap\cap_{i=1}^m\{\lambda_i(n)=l_i,\tau_i(n)=k_i,  \vec S_{\oo K_{i-1},K_i}=\vec x_i\}) 
\\     \nonumber %
=\psi_0q^n{\mathbf P} ( \oo\nu_0=L_0,\alpha(\nu_0)=K_0,\tilde  S_{K_0}=\vec y) 
\cdot\prod_{i=1}^m {\mathbf P} (\oo\lambda_i=l_i,\oo\tau_i=k_i,
\widetilde S_{K_{i-1}K_i}=\vec x_i\}). 
\end{gather*}
\end{corollary}

Comment that random vectors $\{\oo\xi_j,\ j=1,2\dots\}$ that were introduced in (\ref{b18}) may be dependent, notwithstanding that $\{\xi_j,\ j=1,2\dots\}$ were i.i.d.
However, the random blocks 
$$\big(\oo{\lambda}_i, \oo \tau_i, (\oo{\xi}_{K_{i-1}+1},\ldots,\oo{\xi}_{K_i}) \big),\ i=1,2\dots,$$
 are i.i.d. and do not depend on the initial block $\left( \overline{S}_0, \oo{\nu}_0, \oo T_0,\, (\oo\xi_1,\ldots, \oo\xi_{\oo{T}_0}) \right)$.
This type of the phenomenon is typical for conditioning that involves infinite future: an i.i.d. sequence is transformed into a regenerative sequence.
It appears even in the simplest scenario, for a one-dimensional random walk with positive drift, 
conditioned to stay positive (see,e.g.,  \cite{FZ},\cite{Ku},\cite{MS} for similar observations in ``unconditioned'' models).

\end{document}